\documentclass[11pt]{article}
\usepackage{graphicx}

\newfont{\bb}{msbm10 at 11pt}
\def\r{\hbox{\bb R}}

\newcommand{\s}{\hbox{\bf S}}
\newcommand{\e}{\hbox{\bf E}}
\newcommand{\h}{\hbox{\bf H}}
\newtheorem{theorem}{Theorem}[section]

\newtheorem{definition}{Definition}[section]

\begin{document}

\title{Non-degenerate surfaces of revolution  in Minkowski space that satisfy the relation $aH+bK=c$}
\author{\"{O}zg\"{u}r Boyac\i o\u{g}lu Kalkan \\
Mathematics Department\\
Afyon Kocatepe University\\
Afyon 03200 Turkey\\
email: bozgur@aku.edu.tr\\
 \vspace*{.5cm}\\
Rafael L\'opez\thanks{%
Partially supported by MEC-FEDER grant no. MTM2007-61775 and Junta de Andaluc\'{\i}a grant no. P06-FQM-01642.}\\
Departamento de Geometr\'{\i}a y Topolog\'{\i}a\\
Universidad de Granada\\
18071 Granada, Spain\\
email: rcamino@ugr.es\\
 \vspace*{.5cm}\\
Derya Saglam\\ Mathematics Department\\
 Afyon Kocatepe Universtiy\\
  Afyon 03200 Turkey}
\date{}
\maketitle

\begin{abstract}
In this work, we study spacelike and timelike surfaces of revolution in Minkowski space $\e_{1}^{3}$  that satisfy $aH+bK=c$, where $H$ and $K$ denote the mean curvature and the Gauss curvature of the surface and $a$, $b$ and $c$ are constants. The classification depends on the causal character of the axis of revolution and in all the cases, we obtain a first integral of the equation of the generating curve of the surface.
\end{abstract}


\section{Introduction}\label{intro}

Consider the three-dimensional Minkowski space $\e_1^3$, that is, the real vector space $\r^3$ endowed with the Lorentzian metric $\langle,\rangle=(dx)^2+(dy)^2-(dz)^2$, where $(x,y,z)$ stand for the usual coordinates of $\r^3$.
A vector $v\in\e_1^3$ is said spacelike if $\langle v,v\rangle>0$
or $v=0$, timelike if $\langle v,v\rangle<0$ and lightlike if $\langle
v,v\rangle=0$ and $v\not=0$. A submanifold $S\subset\e_1^3$  is said
spacelike, timelike or lightlike if the induced metric on $S$ is a
Riemannian metric (positive definite), a Lorentzian metric (a metric of
index $1$) or a degenerated metric, respectively. In the case that $S$ is a straight-line $L=<v>$, this means that $v$ is spacelike, timelike or lightlike, respectively. If $S$ is a plane $P$, this is equivalent that
any orthogonal vector to $P$ is timelike, spacelike or lightlike
respectively.

An immersion $x:M\rightarrow\e_1^3$  of a surface $M$ is said non-degenerated  if the induced metric $x^*(\langle,\rangle)$ on $M$ is non-degenerate. In this setting, there is only two possibilities: if $x^*(\langle,\rangle)$ is positive definite , that is, it is a Riemmannian metric and the immersion is called {\it spacelike} or $x^*(\langle,\rangle)$ is   a Lorentzian metric, that is, a metric of index $1$, and the immersion is called {\it timelike}. For spacelike surfaces, the tangent planes are spacelike everywhere, and for timelike surfaces, they are timelike.

We consider spacelike or timelike surfaces in $\e_1^3$ that satisfy the relation
\begin{equation}\label{w1}
aH+bK=c,
\end{equation}
where $H$ and $K$ are the mean curvature and the Gauss curvature of the surface, and $a$, $b$ and $c$ are constants. We say that the surface is a {\it linear Weingarten surface} of $\e_1^3$. In general, a Weingarten surface is a surface that satisfies a certain smooth relation $W=W(H,K)=0$ and our case, that is, surfaces that satisfy (\ref{w1}) is the simplest case of $W$, that is, that $W$ is a linear function in its variables. The family of linear Weingarten surfaces include the surfaces with constant mean curvature ($b=0$) and the surfaces with constant Gauss curvature ($a=0$).

In this work we study linear Weingarten surfaces that are rotational, that is, invariant by a group of motions of $\e_1^3$ that pointwised fixed a straight-line. In such case, Equation (\ref{w1}) is a second ordinary differential equation that describes the shape of the generating curve of the surface. One can not expect to integrate this equation, because even in the trivial cases that $a=0$ or $b=0$, this integration is not possible.  We are going to discard the cases that $H$ is constant of $K$ is constant, which are known: see for example \cite{hn,lo1,lo2}. We will obtain a first integration of (\ref{w1}). For the particular case that $a^2-4bc\epsilon=0$,  we describe all solutions, exactly, we have

\begin{theorem} Let $M$ be a non-degenerate rotational surface in $\e_1^3$, and take $\epsilon=1$ if $M$ is spacelike and $\epsilon=-1$ if $M$ is timelike. Assume that $M$ is a linear Weingarten surface such that $a^2-4bc\epsilon=0$. After a rigid motion of the ambient space, a parametrization $X(u,v)$ of $M$ is  as follows:
\begin{enumerate}
\item If the axis if timelike, $X(u,v)=(u\cos(v),u\sin(v),z(u))$, where
$$z(u)=\pm  \sqrt{\frac{4\epsilon b^2}{a^2}+(\frac{C}{a}\pm u)^2}+\mu,\ \ C=2\sqrt{ b\epsilon(-b+\lambda)},\  \mu,\lambda\in\r.$$
\item If the axis is spacelike, we have two possibilities:
\begin{enumerate}
\item The parametrization is $X(u,v)=(u,z(u)\sinh(v),z(u)\cosh(v))$, where
$$z(u)=  \pm \frac{C}{a}\pm\sqrt{\frac{4\epsilon b^2}{a^2}\pm (u\pm \mu)^2},\  \ C=2\sqrt{ b\epsilon \lambda},\  \mu,\lambda\in\r.$$
\item The parametrization is $X(u,v)=(u,z(u)\cosh(v),z(u)\sinh(v))$, where
$$z(u)=    \frac{-C}{a}\pm\sqrt{\frac{4  b^2}{a^2}\pm (u\pm \mu)^2},\ \ C=2\sqrt{ b \lambda},\  \mu,\lambda\in\r.$$
\end{enumerate}
\item If the axis is lightlike, $X(u,v)=(-2 uv,z(u)+u -u v^2,z(u)-u-u v^2)$, where
$$z(u)=\frac{1}{48}\Big(\frac{-4ac\lambda+(cC^2-2a^2\lambda)u}{\epsilon c\lambda(2\lambda+cu^2)}+\epsilon
\frac{cC^2+2a^2\lambda}{\sqrt{-2c\lambda}}\mathrm{arc}\tanh(\sqrt{-\frac{c}{2\lambda}}\ u)\Big)+\mu,\ \  \mu,\lambda\in\r.$$
\end{enumerate}
\end{theorem}


\section{Rotational  surfaces in $\e_1^3$}


In this section we  describe the surfaces of revolution of $\e_1^3$ and we recall the concepts of mean curvature and Gauss curvature for a non-degenerate surface. We consider the rigid motions of the ambient space that leave
a straight-line pointwised fixed, called, the axis of the surface.  Let $L$ be the axis of the surface. Depending on $L$, there are three types of rotational motions. After an isometry of $\e_1^3$, the expressions of rotational motions with respect to the canonical basis  $\{e_1,e_2,e_3\}$  are as follows:
$$R_v:\left(\begin{array}{c}x_1\\ x_2\\ x_3\end{array}\right)\longmapsto
\left(\begin{array}{ccc}\cos{v}&\sin{v}&0\\-\sin{v}&\cos{v}&0\\0&0&1\end{array}\right)\left(\begin{array}{c}x_1\\ x_2\\ x_3\end{array}\right).$$
$$R_v: \left(\begin{array}{c}x_1\\ x_2\\  x_3\end{array}\right)\longmapsto
\left(\begin{array}{ccc}1&0&0\\0&\cosh{v}&\sinh{v}\\0&\sinh{v}&\cosh{v}\end{array}\right)\left(\begin{array}{c}x_1\\ x_2\\ x_3\end{array}\right).$$
$$R_v: \left(\begin{array}{c}x_1\\ x_2\\ x_3\end{array}\right)\longmapsto
\left(\begin{array}{ccc}1&-v&v\\ v&1-\frac{v^2}{2}&\frac{v^2}{2}\\ v&-\frac{v^2}{2}&1+\frac{v^2}{2}
\end{array}\right)\left(\begin{array}{c}x_1\\ x_2\\ x_3\end{array}\right).$$
See \cite{on,we} for more details.

\begin{definition} A surface $M$ in $\e_1^3$ is a surface of revolution, or
rotational surface, if  $M$ is
invariant by some of the above three groups of rigid motions.
\end{definition}
In particular, there exists a planar curve $\alpha=\alpha(u)$ that generates the surface, that is, $M$ is the set of points given by
$\{R_v(\alpha(u));u\in I,v\in\r\}$. We now describe the parametrizations of a  rotational surface.
\begin{enumerate}
\item \emph{Case $L$ is a timelike axis}. Consider that  $L$ is the $x_3$-axis. If $p=(x_0,y_0,z_0)\not\in L$, then $\{R_v(p);v\in\r\}$ is an Euclidean circle of radius
$\sqrt{x_0^2+y_0^2}$ in the plane $z=z_0$. If $\alpha(u)=(u,0,z(u))$ is a planar curve in the plane $y=0$, then the surface of revolution generated by $\alpha$ writes as
\begin{equation}\label{rot1}
X(u,v)=(u\cos(v),u\sin(v),z(u)),\  u\not=0.
\end{equation}
\item \emph{Case $L$ is a spacelike axis}.  Consider that $L$ is the $x_1$-axis. If $p=(x_0,y_0,z_0)$ does not belong to $L$, then $\{R_v(p);v\in\r\}$ is an Euclidean hyperbola in the
plane $x=x_0$ and with equation $y^2-z^2=y_0^2-z_0^2$.  For this kind of rotational surfaces, we have two type of surfaces:
\begin{enumerate}
\item If $\alpha(u)=(u,0,z(u))$ is a planar curve in the plane $y=0$, then the surface of revolution generated by $\alpha$ writes as
\begin{equation}\label{rot21}
X(u,v)=(u,z(u)\sinh(v),z(u)\cosh(v)),\   u\not=0.
\end{equation}
\item If $\alpha(u)=(u,z(u),0)$ is a planar curve in the plane $z=0$, then the surface is given by
\begin{equation}\label{rot22}
X(u,v)=(u,z(u)\cosh(v),z(u)\sinh(v)),\   u\not=0.
\end{equation}
\end{enumerate}
\item \emph{Case $L$ is a lightlike axis}. Consider that $L$ is the  straight-line $v_1=<(0,1,1)>$. If $p=(x_0,y_0,z_0)$ does not belong to the plane $<e_1,v_1>$,    the orbit $\{R_v(p);v\in\r\}$ is the curve
 $$\beta(v)=(x-(y-z)v,xv+y-(y-z)\frac{v^2}{2},xv+z-(y-z)\frac{v^2}{2}).$$
 The curve $\beta$ lies in the plane $y-z=y_0-z_0$ and describes a parabola in this plane, namely,
$$\beta(v)=(x,y,z)+v(-(y-z) e_1+x v_1)-\frac{y-z}{2}v^2v_1.$$
Consider $\alpha(u)$ a planar curve in the plane $<(0,1,1),(0,1,-1)>$ given  as a graph on
the straight-line $<(0,1,-1)>$, that is, $\alpha(u)=(0,u+z(u),-u+z(u))$. The surface of revolution generated by $\alpha$ is
\begin{equation}\label{rot3}
X(u,v)=(-2 uv,z(u)+u -u v^2,z(u)-u-u v^2),\  u\not=0.
\end{equation}
\end{enumerate}

Let $M$ be   surface and $x:M\rightarrow \e_1^3$ a non-degenerate immersion and we simply say that $M$ is non-degenerate. The surface could be not orientable, but if the immersion is spacelike, then $M$ is necessarily orientable. This is due to the following fact. At each point $p\in M$ there is two possible choices of a unit normal vector to the tangent plane $T_pM$ of $M$ at $p$. The normal vector to $M$ is a timelike vector, and in Minkowski space, two any timelike vectors are not orthogonal. Thus, if $E_3=(0,0,1)$, at each point $p\in M$, we take that unit normal vector $N(p)$ such that $\langle N(p),E_3\rangle<0$. This allows to define an global orientation on $M$, proving that $M$ is orientable. With this choice of $N$, we say that $N$ is future directed.   In the case that the immersion is timelike, we will assume that $M$ is orientable.

Let $x:M\rightarrow\e_1^3$ be a non-degenerate immersion of a surface  $M$ and let  $N$ be a Gauss map.
Let $U, V$ be  vector fields to $M$ and we denote by  $\nabla^0$ and $\nabla$ the Levi-Civitta connections of  $\e_1^3$ and $M$ respectively.  The Gauss formula says
$\nabla_U^0 V=\nabla_U V+\mbox{II}(U,V),$
where $\mbox{II}$ is the second fundamental form of the immersion. The Weingarten endomorphism is $A_p:T_pM
\rightarrow T_p M$ defined as  $A_p(U)=-(\nabla_U^0 N)_p^\top=(-dN)_p(U)$.   We have then
 $\mbox{II}(U,V)=-\epsilon\langle \mbox{II}(U,V),N\rangle N=-\epsilon\langle AU,V\rangle N$, where $\epsilon=1$ if $M$ is spacelike and $\epsilon=-1$ if $M$  is timelike.  The mean curvature vector $\vec{H}$ is defined as
$\vec{H}=(1/2)\mbox{trace}(\mbox{II})$ and the Gauss curvature $K$ as the determinant of $\mbox{II}$ computed in both cases with respect to an orthonomal basis. The mean curvature $H$ is the function given by $\vec{H}=HN$, that is, $H=-\epsilon \langle\vec{H},N\rangle$. If $\{e_1,e_2\}$ is an orthonormal basis at each tangent plane, with $\langle e_1,e_1\rangle=1$, $\langle e_2,e_2\rangle=\epsilon$,  then
\begin{eqnarray*}\vec{H}&=&\frac12\mbox (\mbox{II}(e_1,e_1)+\mbox{II}(e_2,e_2))=-\epsilon\frac12(\langle Ae_1,e_1\rangle+\epsilon\langle Ae_2,e_2\rangle)N=-\epsilon(\frac12 \mbox{trace}(A))N\\
K&=&-\epsilon\mbox{det}(A).
\end{eqnarray*}
In this work we need to compute $H$ and $K$ using a parametrization of the surface. Let  $X:D\subset\r^2\rightarrow\e_1^3$ be a parametrization of the surface, $X=X(u,v)$. Then $A=\mbox{II}(\mbox{I})^{-1}$, $\mbox{I}=\langle,\rangle$ and
we have the known formulae (\cite{we}):
\begin{equation}\label{hk}
H=-\epsilon\frac12\frac{eG-2fF+gE}{EG-F^2},\hspace*{1cm}K=-\epsilon\frac{eg-f^2}{EG-F^2},\end{equation}
where $\{E,F,G\}$ and $\{e,f,g\}$ are the coefficients of $\mbox{I}$ and $\mbox{II}$, respectively:
$$E=\langle X_u,X_u\rangle,\ F=\langle X_u,X_v\rangle,\ G=\langle X_v,x_v\rangle,$$
$$e=-\langle N_u,X_{u}\rangle,\ f=-\langle N_u,X_{v}\rangle,\ g=-\langle N_v,X_{v}\rangle,$$
where the subscripts denote the corresponding derivatives. Here $N$ is
$$N=\frac{X_u\times X_v}{\sqrt{\epsilon(EG-F^2)}}.$$
 We recall that
$$W:=EG-F^2=\epsilon|X_u\times X_v|^2\ \left\{\begin{array}{l}
\mbox{is positive if $M$ is spacelike}\\
\mbox{is negative if $M$ is timelike}\end{array}\right.$$
Finally, in order to the computations for $H$ and $K$, we recall that the cross-product $\times$ satisfies that
for any vectors $u,v,w\in\e_1^3$, $\langle u\times v,w\rangle=\mbox{det}(u,v,w)$. Then
(\ref{hk}) writes as
\begin{equation}\label{hk21}
H=-\frac{\epsilon}{2}\frac{G\mbox{det}(X_u,X_v,X_{uu})-
2F\mbox{det}(X_u,X_v,X_{uv})+E\mbox{det}(X_u,X_v,X_{vv})}{(\epsilon(EG-F^2))^{3/2}}:=\frac{H_1}{2W^{3/2}}.
\end{equation}
\begin{equation}\label{hk22}K=-\frac{\mbox{det}(X_u,X_v,X_{uu})\mbox{det}(X_u,X_v,X_{vv})-
\mbox{det}(X_u,X_v,X_{uv})^2}{ (EG-F^2)^2}:=\frac{K_1}{W^2}.
\end{equation}
In Minkowski ambient space, the role of   spheres is played by pseudohyperbolic surfaces and pseudospheres \cite{on}. If $p_0\in \e_1^3$ and $r>0$ the pseudohyperbolic surface centered at $p_0$ with radius $r>0$ is $\h^{2,1}(r;p_0)=\{p\in\e_1^3;\langle p-p_0,p-p_0\rangle=-r^2 \}$ and the pseudosphere centered at $p_0$ and radius $r>0$ is $\s^{2,1}(r;p_0)=\{p\in\e_1^3;\langle p-p_0,p-p_0\rangle=r^2\}$.
  If $M$ is spacelike (resp. timelike) then $N$ is timelike (resp. spacelike) and $N:M\rightarrow \h^{2,1}(1)$ (resp. $N:M\rightarrow \s^{2,1}(1)$), where $\h^{2,1}(1)=\h^{2,1}(1;O)$ (resp. $\s^{2,1}(1)=\s^{2,1}(r;O)$, being $O$ the origin of coordinates of $\r^3$.  For both kind of surfaces, we can take  $N(p)=(p-p_0)/r$ and $A=-\frac{1}{r}I$. Then  $H=\epsilon/r$ and $K=-\epsilon/r^2$.
\section{Rotational surfaces with timelike axis}

We assume that the generating curve $\alpha$ lies in the $xz$-plane and we parametrize $\alpha$ as the graph of a function $z=z(u)$, that is,   $\alpha(u)=(u,0,z(u))$, $u>0$. Then the surface is parametrized as in (\ref{rot1}) and
$W=u^2(1-z^2)$. Thus $z'^2<1$ if the surface is spacelike and $z'^2>1$ if $M$ is timelike. Using (\ref{hk21}) and (\ref{hk22}), the expressions of $H$ and $K$ are:
$$H=-\frac12\Bigg(\frac{\epsilon z'}{u\sqrt{\epsilon(1-z'^2)}}+\frac{z''}{(\epsilon(1-z'^2))^{3/2}}\Bigg),\ \ K=-\frac{z'z''}{ u(1-z'^2)^2}.$$
Then the relation (\ref{w1}) writes as
$$\frac{a}{2}\Bigg(\frac{\epsilon z'}{u\sqrt{\epsilon(1-z'^2)}}+\frac{z''}{(\epsilon(1-z'^2))^{3/2}}\Bigg) +b\frac{z'z''}{ u(1-z'^2)^2}=-c.$$
Multiplying by $u$ we obtain a first integral. Exactly, we have
$$  a \Bigg( u \frac{\epsilon z'}{\sqrt{\epsilon(1-z'^2)}}\Bigg)'+b\Bigg(\frac{1}{1-z'^2}\Bigg)'=-2cu.$$
Then there exists a integration constant $\lambda\in\r$ such that
\begin{equation}\label{int1}
\epsilon  \frac{a uz'}{\sqrt{\epsilon(1-z'^2)}}+ \frac{b}{ 1-z'^2}=-c u^2+\lambda.
\end{equation}
Let
$$\phi=\frac{z'}{\sqrt{\epsilon(1-z'^2)}}.$$
Then $1+\epsilon\phi^2=1/(1-z'^2)$ and Equation (\ref{int1}) writes as
$b\phi^2+au\phi+\epsilon (b+c u^2-\lambda)=0$. Hence, we obtain $\phi$:
\begin{equation}\label{sol1}
\frac{z'}{\sqrt{\epsilon(1-z'^2)}}= \frac{-au\pm\sqrt{(a^2-4bc\epsilon )u^2+4b\epsilon(-b+\lambda)}}{2b}.
\end{equation}
We completely solve this differential equation in two particular cases:
\begin{enumerate}
\item Consider $\lambda= b$. Then we have
$$\frac{z'}{\sqrt{\epsilon(1-z'^2)}}= \frac{-a \pm\sqrt{a^2-4bc\epsilon   }}{2b}u= Cu,\ \ C=\frac{-a \pm\sqrt{a^2-4bc\epsilon   }}{2b}.$$
Then
$$z(u)=\pm\frac{\sqrt{\epsilon+C^2 u^2}}{C}+\mu,\ \mu\in\r.$$ From the parametrization (\ref{rot1}) of the surface, one concludes that  $M$ satisfies the equation $x^2+y^2-(z-\mu)^2=-\frac{\epsilon}{C^2}$. Letting
$p_0=(0,0,\mu)$,  if $\epsilon=1$, the surface  $M$ is the pseudohyperbolic surface $\h^{2,1}(1/|C|;p_0)$ and when $\epsilon =-1$, $M$ is a pseudosphere $\s^{2,1}(1/|C|;p_0)$.
\item Assume $a^2-4bc\epsilon =0$. Then
$$\frac{z'}{\sqrt{\epsilon(1-z'^2)}}= \frac{-au\pm C}{2b},\ \ C=2\sqrt{ b\epsilon(-b+\lambda)}.$$
The integration of this equation is
$$z(u)=\pm  \sqrt{\frac{4\epsilon b^2}{a^2}+(\frac{C}{a}\pm u)^2} +\mu,\ \mu\in\r.$$

\end{enumerate}

\section{Rotational surfaces with spacelike axis}

We distinguish two cases according the two possible parametrizations.

\begin{enumerate}
\item Case I. Assume that the parametrization is given by (\ref{rot21}). The relation (\ref{w1}) writes as
$$\frac{a}{2} \Bigg(\frac{\epsilon}{z\sqrt{\epsilon(1-z'^2)}}+\frac{z''}{(\epsilon(1-z'^2))^{3/2}}\Bigg)+b\frac{z''}{z(1-z'^2)^2}=-c.$$
Multiplying by $zz'$, we obtain a first integral. Exactly, we have
$$a   \Bigg( \frac{\epsilon z}{\sqrt{\epsilon(1-z'^2)}}\Bigg)'+b\Bigg(\frac{1}{1-z'^2}\Bigg)'=-c(z^2)'.$$
Then there exists an integration constant $\lambda\in\r$ such that
\begin{equation}\label{int2}
\epsilon  \frac{a z}{\sqrt{\epsilon(1-z'^2)}}+ \frac{b}{1-z'^2}=-cz^2+\lambda.
\end{equation}
Now we take $\phi=1/\sqrt{\epsilon(1-z'^2)}$. Then Equation (\ref{int2}) writes as
$$b\phi^2+a z\phi+ \epsilon(cz^2-\lambda)=0.$$
Then
\begin{equation}\label{sol2}
\frac{1}{\sqrt{\epsilon(1-z'^2)}}= \frac{-a  z\pm\sqrt{(a^2-4bc\epsilon)z^2+4b\epsilon\lambda}}{2b}.
\end{equation}
We completely solve this differential equation in two particular cases:
\begin{enumerate}
\item Consider $\lambda= 0$. Then we have
$$\frac{1}{\sqrt{\epsilon(1-z'^2)}}= \frac{-a  \pm\sqrt{a^2-4bc\epsilon}}{2b}z= Cz,\ \ C=\frac{-a  \pm\sqrt{a^2-4bc\epsilon}}{2b}.$$
The solution of this differential equation is
$$z(u)=\pm \sqrt{\frac{\epsilon}{C^2}\pm (u\pm C\mu)^2},\ \mu\in\r\}.$$ From the parametrization (\ref{rot21}) of the surface, one concludes that $M$ satisfies the equation $(x-C\mu)^2+y^2-z^2=-\frac{\epsilon}{C^2}$. Thus, if we set $p_0=(\pm C\mu,0,0)$, for $\epsilon=1$ we obtain that  $M$ is the pseudohyperbolic surface
$\h^{2,1}(1/|C|;p_0)$ and for $\epsilon =-1$, $M$ is the pseudosphere $\s^{2,1}(1/|C|;p_0)$.
\item Assume $a^2-4bc\epsilon =0$. Then
$$\frac{1}{\sqrt{\epsilon(1-z'^2)}}= \frac{-az\pm C}{2b},\ \ C=2\sqrt{ b\epsilon \lambda}.$$
The integration of this equation is
$$z(u)=  \pm \frac{C}{a}\pm\sqrt{\frac{4\epsilon b^2}{a^2}\pm (u\pm \mu)^2},\ \mu\in\r.$$
\end{enumerate}
 \item Case II. The expression of the parametrization is written in (\ref{rot22}). In this case, the surface is timelike, since $EG-F^2=-z^2(1+z'^2)$. The Weingarten relation (\ref{w1}) is
$$\frac{a}{2} \Bigg(\frac{-1}{z\sqrt{ 1+z'^2}}+\frac{z''}{( 1+z'^2)^{3/2}}\Bigg)-b\frac{z''}{z(1+z'^2)^2}=c.$$
Multiplying by $zz'$ again, we have
$$-a   \Bigg( \frac{  z}{\sqrt{ 1+z'^2}}\Bigg)'+b\Bigg(\frac{1}{1+z'^2}\Bigg)'=c(z^2)'.$$
It follows the existence of an integration constant $\lambda\in\r$ such that
\begin{equation}\label{int22}
-  \frac{a z}{\sqrt{ 1+z'^2}}+ \frac{b}{1+z'^2}=cz^2+\lambda.
\end{equation}
If we set $\phi=1/\sqrt{ 1+z'^2}$,  Equation (\ref{int22}) is
$b\phi^2-a z\phi -cz^2-\lambda=0$, obtaining
\begin{equation}\label{sol22}
\frac{1}{ 1+z'^2}= \frac{a  z\pm\sqrt{(a^2+4bc)z^2+4b\lambda}}{2b}.
\end{equation}
As in the previous case, we solve this equation in the next two  cases:
\begin{enumerate}
\item If $\lambda= 0$, then
$$\frac{1}{\sqrt{ 1+z'^2}}= \frac{-a  \pm\sqrt{a^2+4bc}}{2b}z= Cz,\ \ C=\frac{a  \pm\sqrt{a^2+ 4bc}}{2b}.$$
The solution of this  equation is
$$z(u)=\pm \sqrt{\frac{1}{C^2}-(u\pm C\mu)^2},\ \mu\in\r\}.$$
This surface is the pseudosphere  $\s^{2,1}(1/|C|;p_0)$, with $p_0=(\pm C\mu,0,0)$ since by the expression of the parametrization (\ref{rot22}), the coordinates of $M$ satisfies $(x\pm C\mu)^2+y^2-z^2=1/C^2$.
\item If $a^2+4bc =0$, then
$$\frac{1}{\sqrt{1+z'^2}}= \frac{az\pm C}{2b},\ \ C=2\sqrt{ b  \lambda}.$$
The solution of this equation is
$$z(u)=    \frac{-C}{a}\pm\sqrt{\frac{4  b^2}{a^2}\pm (u\pm \mu)^2},\ \mu\in\r.$$
 \end{enumerate}

\end{enumerate}

\section{Rotational surfaces with lightlike axis}
Consider the parametrization given in (\ref{rot3}). Then $EG-F^2=16u^2z'$ and the  relation (\ref{w1}) writes as
$$\frac{a}{2}\Bigg(\frac{1}{2u\sqrt{\epsilon z'}}-\frac{\epsilon z''}{4(\epsilon z')^{3/2}}\Bigg)+b\frac{z''}{8uz'^2}=c.$$
Multiplying by $u$ we obtain a first integral. Exactly, we have
$$\frac{a}{4} \Bigg( \frac{  u}{\sqrt{\epsilon z'}}\Bigg)'- \frac{b}{8}\Bigg(\frac{1}{z'}\Bigg)'=cu.$$
Then there exists a integration constant $\lambda\in\r$ such that
\begin{equation}\label{int3}
\frac{a}{4}  \frac{ u}{\sqrt{\epsilon z'}}- \frac{b}{8z'}=\frac{c}{2}u^2+\lambda.
\end{equation}
From (\ref{int3}), we obtain the value of $\sqrt{\epsilon z'}$:
$$\sqrt{\epsilon z'}=\frac{a\epsilon u\pm\sqrt{(a^2-4bc\epsilon)u^2-8b\epsilon\lambda}}{4\epsilon(c u^2+2\lambda)}.$$
As in the two previous cases, we distinguish two special cases:
\begin{enumerate}
\item If $\lambda=0$, then
$$\sqrt{\epsilon z'}=\frac{a\pm\epsilon\sqrt{a^2-4bc\epsilon}}{4c}\frac{1}{u}:=\frac{C}{u},\ \ C=\frac{a\pm\epsilon \sqrt{a^2-4bc\epsilon}}{4c }.$$
We solve this equation obtaining
$$z(u)=-\frac{\epsilon C^2}{u}+\mu,\ \ \mu\in\r.$$
From the parametrization (\ref{rot3}), we see that $M$ satisfies the equation $x^2+y^2-(z-\mu)^2=-4\epsilon C^2$. Thus, if $p_0=(0,0,\mu)$, we have that
$M=\h^{2,1}(2|C|;p_0)$ if $\epsilon=1$, and $M=\s^{2,1}(2|C|;p_0)$ if $\epsilon=-1$.
\item Assume $a^2-4bc\epsilon=0$. Then
$$\sqrt{\epsilon z'}=\frac{a\epsilon u\pm C}{4\epsilon  (cu^2+2\lambda)},\ \ C=\sqrt{-8b\epsilon\lambda}.$$
We point out that $-8b\epsilon\lambda>0$ and that combining with   $a^2-4bc\epsilon=0$, we have $c\lambda\leq 0$.
The solution is
$$z(u)=\frac{1}{64}\Big(\frac{\mp 4aC\lambda\pm\epsilon(cC^2-2a^2\lambda)u}{\epsilon c\lambda(2\lambda+cu^2)}+\epsilon
\frac{cC^2+2a^2\lambda}{\sqrt{-2c^3\lambda^3}}\mathrm{arc}\tanh(\sqrt{-\frac{c}{2\lambda}}\ u)\Big)+\mu,\ \  \mu,\lambda\in\r.$$

\end{enumerate}

\begin{figure}[hbtp]
\begin{center}\includegraphics[width=7cm]{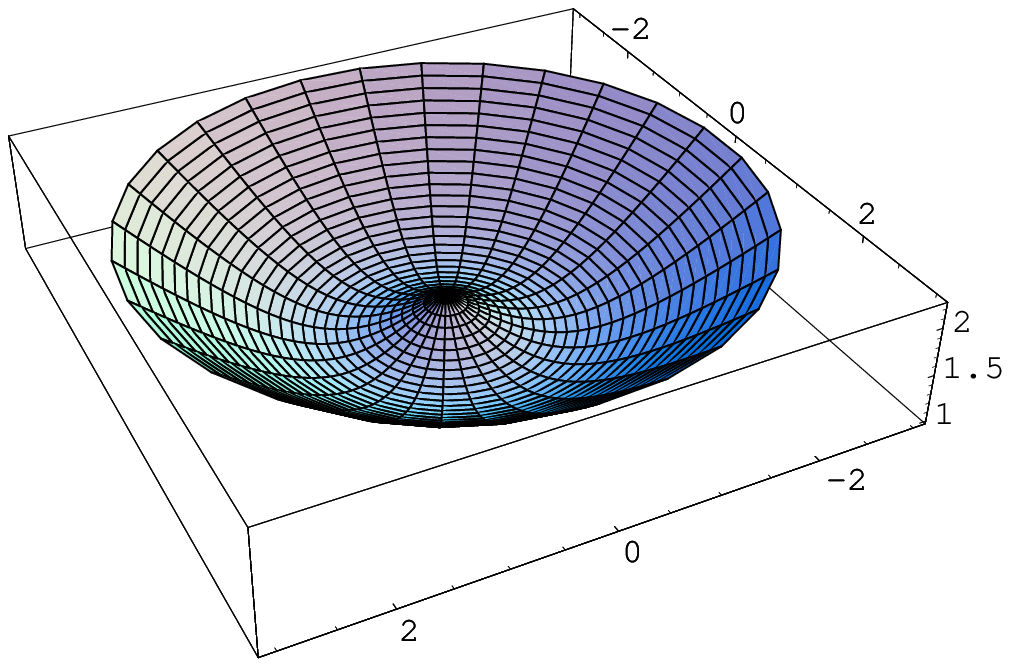}\hspace*{1cm}\includegraphics[width=5cm]{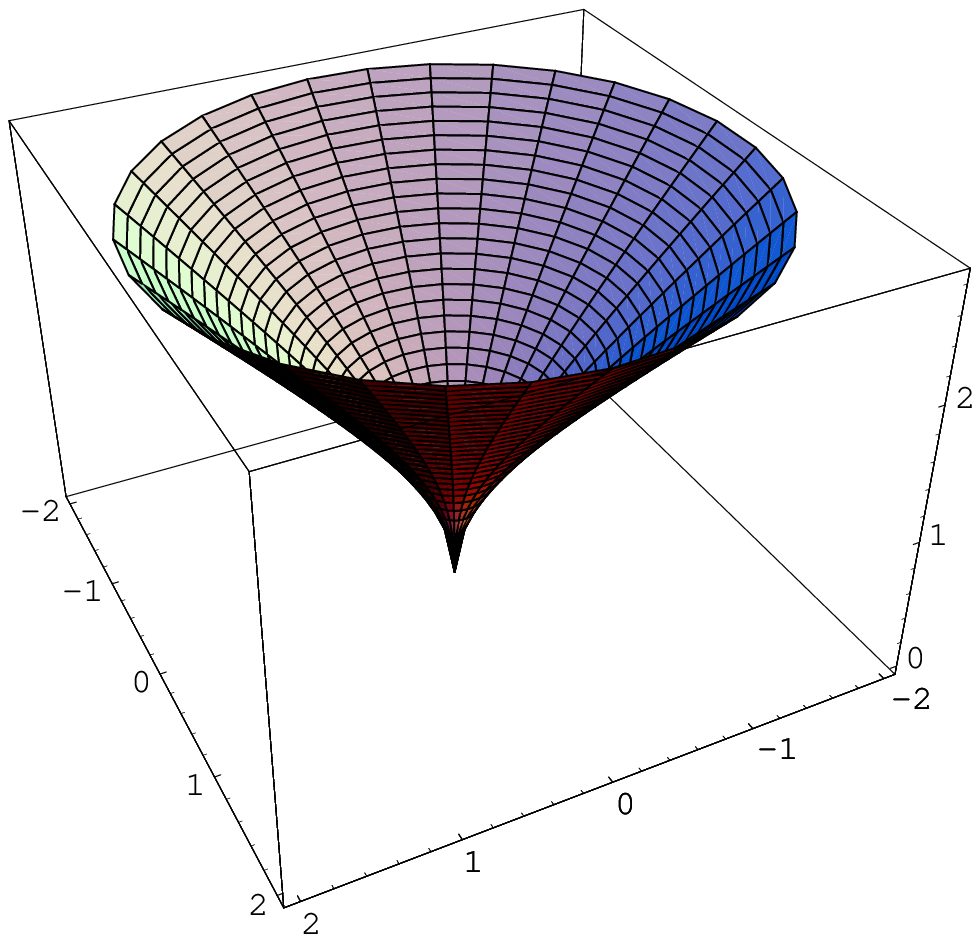}\end{center}
\caption{Rotational surfaces with timelike axis, for $a=2$, $b=\epsilon$ and $\mu=0$: The surface is spacelike with $\lambda=2$ (left). The surface is timelike with $\lambda=0$ (right).}
\end{figure}

\begin{figure}[hbtp]
\begin{center}\includegraphics[width=5cm]{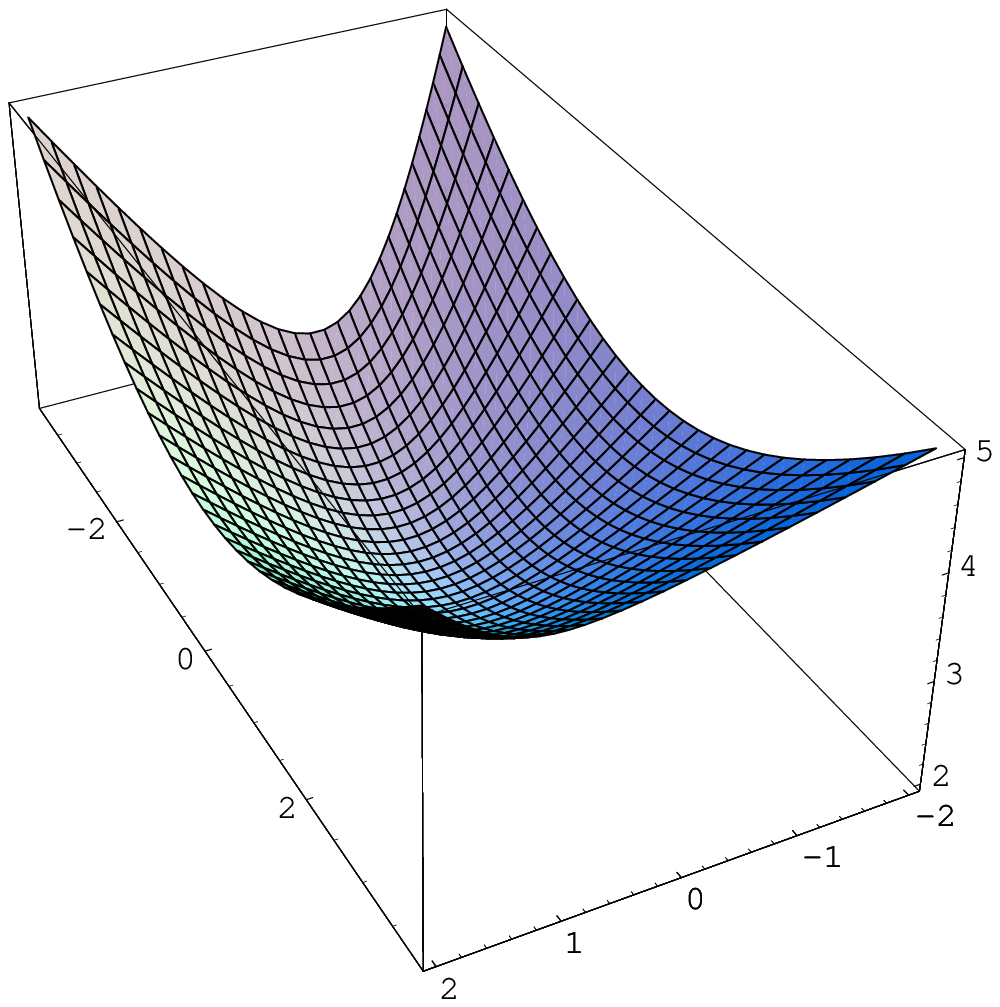}\hspace*{1cm}\includegraphics[width=5cm]{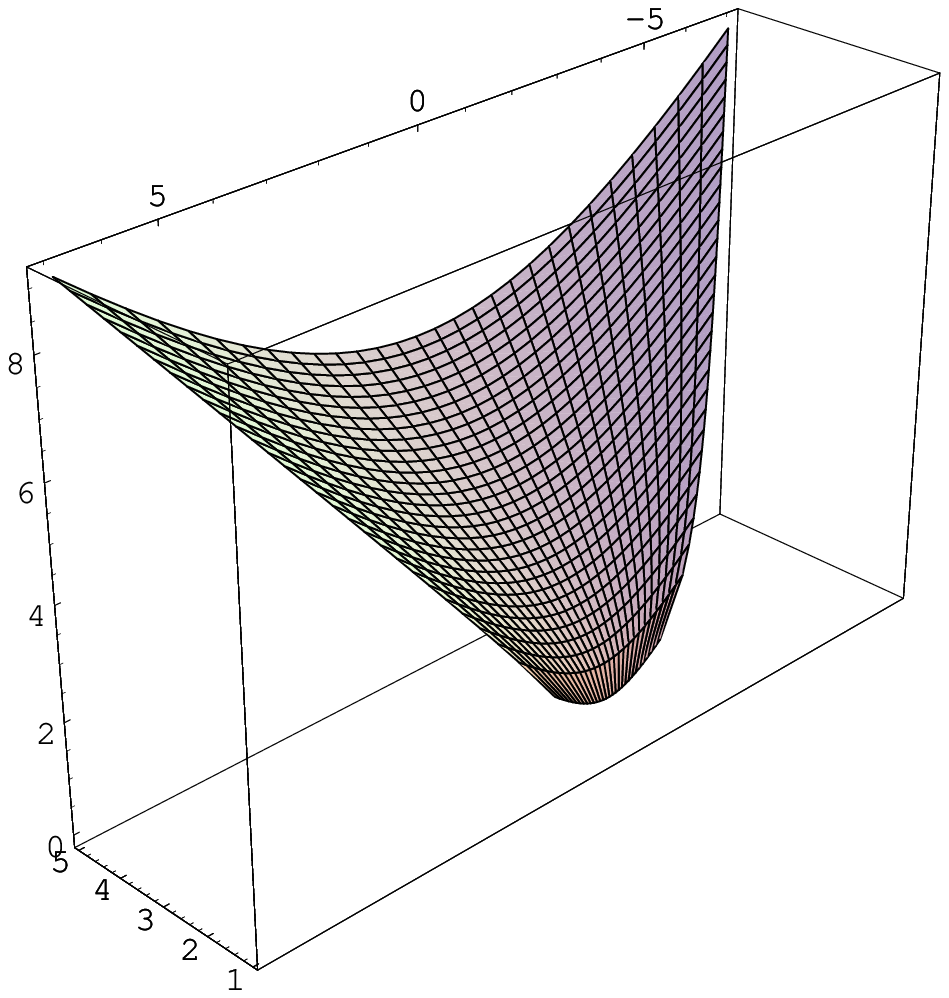}\end{center}
\caption{Rotational surfaces with spacelike axis, for $a=2$, $b=\epsilon$, $\lambda=1$ and $\mu=0$:  The surface is spacelike (left). The surface is timelike (right). }
\end{figure}

\begin{figure}[hbtp]
\begin{center}\includegraphics[width=7cm]{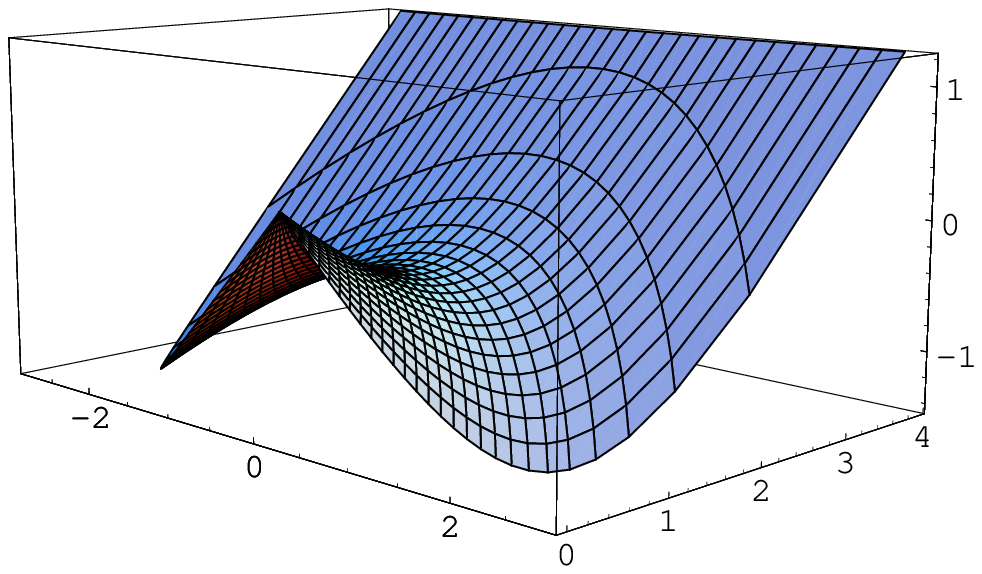}\hspace*{1cm}\includegraphics[width=5cm]{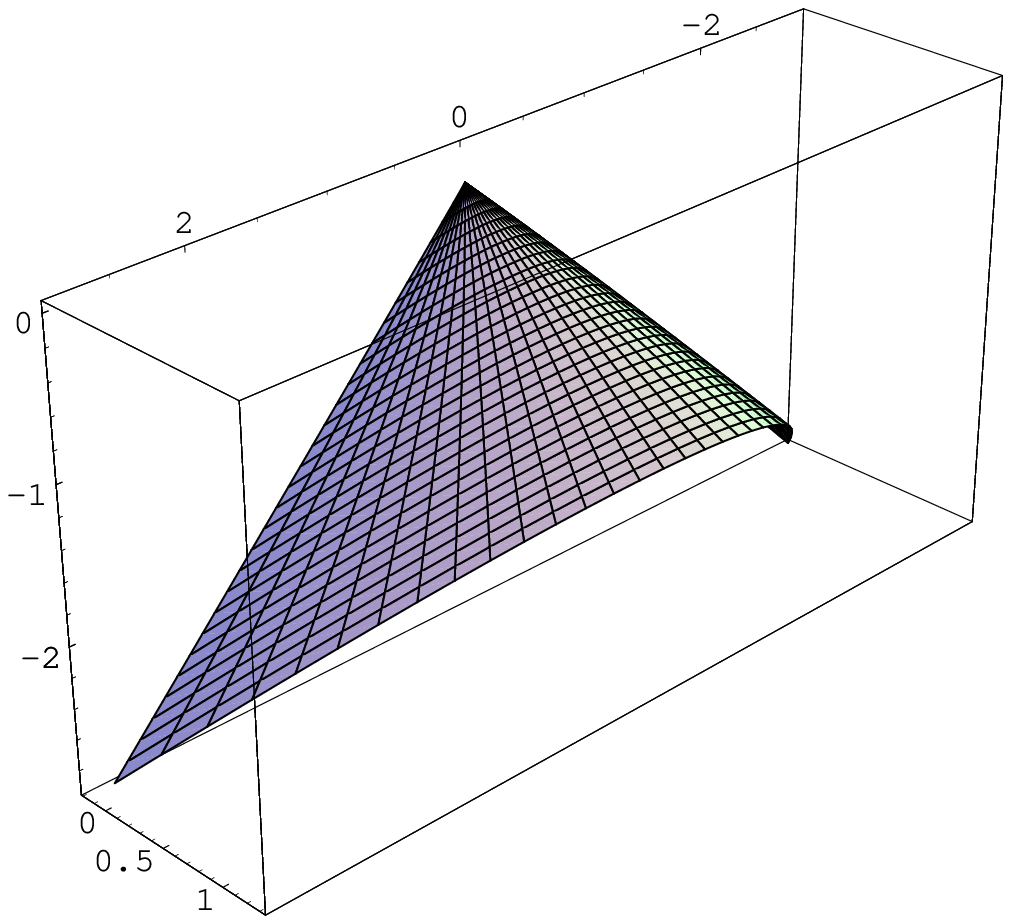}\end{center}
\caption{Rotational surfaces with lightlike axis, for $a=2$, $b=-\epsilon$, $\lambda=1$ and $\mu=0$:  The surface is spacelike (left). The surface is timelike (right). }
\end{figure}


\begin{thebibliography}{99}
\bibitem{hn}  J.I. Hano, K. Nomizu, Surfaces of revolution with constant mean curvature in Lorentz-
Minkowski space, Tohoku Math. J., 36 (1984), 427–-435.

\bibitem{lo1} R.  L\'opez,
Timelike surfaces in Lorentz 3-space with constant mean curvature,
Tohoku Math. J., 52 (2000), 515--532.

\bibitem{lo2} R. L\'opez,
Surfaces of constant Gauss curvature in Lorentz-Minkowski 3-space,
Rocky Mount. J.  Math., 33 (2003), 971--993.

\bibitem{on} B. O'Neill, Semi-Riemannian geometry with applications to relativity, Academic
Press, New York, 1983.

\bibitem{we} T. Weinstein, An introduction to Lorentz surfaces, Walter de Gruyter, Berlin,
1995.
\end{thebibliography}
\end{document}